\documentstyle[12pt]{amsart}
\newtheorem{Theorem}{Theorem}[section] 
\newtheorem{Lemma}[Theorem]{Lemma}
\newtheorem{Corollary}[Theorem]{Corollary}
\newtheorem{Proposition}[Theorem]{Proposition}
\newtheorem{Example}[Theorem]{Example}

\def\Ass{\operatorname{Ass}}
\def\td{\operatorname{tr.deg}}
\def\bideg{\operatorname{bideg}}
\def\Proj{\operatorname{Proj}} 
\def\rdim{\operatorname{rdim}}
\def\sdim{\operatorname{{\it x}-dim}}

\def\To{\longrightarrow} 
\def\mm{{\frak m}} 
\def\PP{{\Bbb P}} 
\def\NN{{\Bbb N}}
\def\ZZ{{\Bbb Z}} 

\begin{document}

\title{Hilbert polynomials\\ of non-standard bigraded algebras}
\author{Nguy\^en Duc Hoang}
\address{Department of Mathematics, Hanoi Pedagogical University, Vietnam}
\author{Ng\^o Vi\^et Trung}
\address{Institute of Mathematics,  Box 631, B\`o H\^o, 10000 Hanoi, Vietnam}
\email{nvtrung@@thevinh.ncst.ac.vn}
\thanks{The authors are partially supported by the National Basic Research Program of Vietnam} 
\subjclass{13D40, 13H15}
\keywords{bigraded algebra, Hilbert polynomial, mixed multiplicity, Rees algebra}

\begin{abstract} 
This paper investigates Hilbert polynomials of bigraded algebras which are generated by elements 
of bidegrees $(1,0), (d_1,1),\ldots,(d_r,1)$, where $d_1,\ldots,d_r$ are non-negative integers. 
The obtained results can be applied to study Rees algebras of homogeneous ideals and their diagonal subalgebras.
 \end{abstract} 
\maketitle

\section*{Introduction} 

Let $R = \oplus_{(u,v) \in \NN^2} R_{(u,v)}$ be a finitely generated bigraded algebra over a field $k$. 
The {\it Hilbert function} of $R$ is the function
$$H_R(u,v) := \dim_kR_{(u,v)}.$$
If $R$ is standard bigraded, i.e. $R$ is generated by elements of bidegrees $(1,0)$ and $(0,1)$, 
then $H_R(u,v)$ is equal to a polynomial in $u,v$ for $u,v$ large enough [Ba], [KMV], [W]. 
This fact does not hold if $R$ is not standard bigraded. \par

In this paper we will study the case $R$ is generated by elements of bidegrees $(1,0),(d_1,1),\ldots,(d_r,1)$, 
where $d_1,\ldots,d_r$ are non-negative integers. 
This case was considered first by P.~Roberts in [Ro1] where it is
shown that there exist integers $c$ and $v_0$ such that $H_R(u,v)$ 
is equal to a polynomial $P_R(u,v)$ for $u \ge cv$ and $v \ge v_0$. 
He calls $P_R(u,v)$ the {\it Hilbert polynomial} of the bigraded algebra $R$.
It is worth to notice that Hilbert polynomials of bigraded algebras of 
the above type appear in Gabber's proof of Serre non-negativity conjecture 
(see e.g. [Ro2]) and that the positivity of certain coefficient 
of such a Hilbert polynomial is strongly related to Serre's positivity
conjecture on intersection multiplicities [Ro3]. \par

Using a different method we are able to show that for 
$$d := \max\{d_1,\ldots,d_r\},$$
there exist  integers $u_0, v_0$ such that $H_R(u,v)= P_R(u,v)$ for $u \ge dv + u_0$ 
and $v \ge v_0$ (Theorem \ref{exist}). Furthermore, the total degree and the degree of $P_R(u,v)$ 
in the variable $u$ can be expressed in terms 
of the dimension of certain quotient rings of $R$ (Theorem \ref{degree} and Theorem \ref{partial}). 
These results cover recent results of P.~Roberts in the case $R$ is generated 
by elements of bidegrees  $(1,0),(0,1),(1,1)$ [Ro3].  
\par

Besides the afore mentioned results we also obtain some results on the leading coefficients of the Hilbert polynomial
$P_R(u,v)$. If $P_R(u,v)$ is written in the form
$$P_R(u,v) = \sum_{i=0}^s \frac{e_i(R)}{i!(s-i)!}u^iv^{s-i} + \text{\rm lower-degree terms},$$
where $s = \deg P_R(u,u)$, we can show that $e_i(R)$ is an integer for $i = 0,\ldots,s$ (Theorem \ref{coefficient}). 
The integer $e_i(R)$ can be negative. That is in sharp contrast to the standard bigraded case 
where it is always non-negative [W]. We shall see that if $R$ is a domain or a Cohen-Macaulay ring, 
then $e_i(R) > 0$ for $i = \deg_u P_R(u,v)$ (Proposition \ref{equi}). 
Moreover, we will show that the integer $e_i(R)$ satisfies the associativity formula (Proposition \ref{associative}).

The inspiration for our study also comes from the theory of diagonal 
subalgebras of bigraded Rees algebras [STV], [CHTV]. Let $A = k[X_1,\ldots,X_n]$ be a polynomial ring over a field $k$ 
and $I$ a homogeneous ideal of $A$. One can associate with $I$ the Rees algebra $A[It] := \oplus_{v\ge 0}I^vt^v.$ 
Since $I$ is a homogeneous ideal, $A[It]$ has a natural bigraded structure by setting $A[It] _{(u,v)} = (I^v)_ut^v$ 
for all $(u,v) \in {\Bbb N}^2$. It is obvious that $A[It]$ belongs to the above class of bigraded algebras. 
Let $V$ denote the blow-up of the projective space $\PP_k^{n-1}$ along the subvariety defined by $I$. 
Then $V$ can be embedded into a projective space by the linear system $(I^e)_c$ for any pair of positive integers $e,c$ 
with $c > de$ [CH]. Such embeddings often yield interesting rational varieties such as the Bordiga-White surfaces [Gi],
the Room surfaces [GG] and the Buchsbaum-Eisenbud varieties [GL]. The homogeneous coordinate ring of 
the embedded variety is the subalgebra $k[(I^e)_c]$ of $A$. It has been observed in [STV] and [CHTV] 
that $k[(I^e)_c]$ can be identified as the subalgebra of $A[It]$ along the diagonal $\{(cv,ev)|\ v \in \NN\}$ of $\NN^2$. 
This idea has been employed successfully to study algebraic properties of the embedded variety. 
Since $P_{A[It]}(cv,ev)$ is the Hilbert polynomial of $k[(I^e)_c]$, we may get uniform information 
on all such embeddings of $V$ from $P_{A[It]}(u,v)$. Furthermore, if we fix $v \ge v_0$, 
then $P_{A[It]}(u,v)$ is the Hilbert polynomial of the ideal $I^v$.\par

If $R$ is a standard bigraded algebra, the coefficients $e_i(R)$ of $P_R(u,v)$ are related to the mixed multiplicities 
introduced by Teissier in singularity theory [Te]. Furthermore, one can use general reductions to compute the numbers 
$e_i(R)$ [Re], [Tr2]. If $R$ is not standard bigraded, we do not know any general method for the computation 
of the integers $e_i(R)$. However, if $R$ is the bigraded Rees algebra $A[It]$ of a homogeneous ideal $I$, 
we have some preliminary information on the Hilbert polynomial $P_{A[It]}(u,v)$ (Theorem \ref{Rees}) 
and we will present an effective method for the computation of $H_{A[It]}(u,v)$ (Proposition \ref{approximation}).\par

We will compute the integers $e_i(R)$ explicitly in the following cases:
\begin{itemize}
\item $R$ is a bigraded polynomial ring of the above type (Proposition \ref{polynomial}), 
\item $R$ is the Rees algebra of an ideal generated by a homogeneous regular sequence (Corollary \ref{regular}),
\item $R$ is the Rees algebra of the ideal generated by the maximal minors of a generic $(r-1)\times r$ matrix (Corollary \ref{minor}).
 \end{itemize}

The last two cases follow from a more general result on ideals generated by homogeneous $d$-sequences 
(Theorem \ref{d-sequence}). We would like to point out that  the computation of the Hilbert polynomial 
of the Rees algebra $A[It]$ is usually very difficult (even if $I$ is generated by a regular sequence) 
since it amounts to the computation of the Hilbert polynomials of all powers of $I$. For instance, 
the computation of these polynomials in the case $I$ being generated by a regular sequence of 
two forms in a polynomial ring of three variables was a key result in the study of certain rational surfaces [GGH].  
Our method will yield a simple proof of this result (Example \ref{GGH}).

The paper is organized as follows. In the first three sections we will study the existence, 
the degree and the coefficients of Hilbert polynomials of bigraded algebras, respectively. 
In the fourth section we apply the obtained result to Rees algebras of homogeneous ideals. 
The last section is devoted to the case the given  ideal is generated by a homogeneous $d$-sequence.

For unexplained terminologies and notations we refer the reader to [E]. Unless otherwise specified, all summations will be ranged over non-negative indices.

\section{Existence of Hilbert polynomial} 

The aim of this section is to prove the following result on the existence of Hilbert polynomials of bigraded algebras.

\begin{Theorem} \label{exist} Let $R$ be a finitely generated bigraded algebra over a field $k$. Assume that $R$ is generated by elements of bidegrees  $(1,0),(d_1,1),\ldots,(d_r,1)$, where $d_1,\ldots,d_r$ are non-negative integers. Set $d = \max\{d_1,\ldots,d_r\}.$
There exist integers $u_0,v_0$ such that for  $u\ge dv+u_0$ and $v \ge v_0$,  the Hilbert function $H_R(u,v)$ is equal to a polynomial $P_R(u,v)$.
\end{Theorem}

If $R$ is generated by elements of bidegrees  $(1,0),(0,1),(1,1)$, P.~Roberts has shown that there exist integers $u_0$
and $v_0$ such that the function $\sum_{t=0}^uH_R(u,v)$ is equal to a polynomial for $u \ge v+u_0$ and $v \ge v_0$
[Ro3, Theorem 1]. But this result can be easily deduced from Theorem \ref{exist}.\par

Unlike the standard bigraded case, we do not have $H_R(u,v)=P_R(u,v)$ for $u,v$ large enough. 
In fact, a bigraded algebra as above may have different Hilbert polynomials, 
depending on the range of the variables in $\NN^2$. 

\begin{Example} 
{\rm Let $S = k[X,Y,Z]$ be a bigraded polynomial ring in three indeterminates $X,Y,Z$ with 
$\bideg X = (1,0)$, $\bideg Y = (0,1)$ and $\bideg Z = (1,1)$. It is easy to check that
$$H_S(u,v) = \cases v+1 & \text{if $u \ge v$,} \\ u+1 & \text{if $v \ge u$.} \endcases $$}
\end{Example}

The proof of Theorem \ref{exist} can be reduced to the case of a bigraded polynomial ring. In this case we have a more precise statement as follows.

\begin{Proposition} \label{polynomial} Let $S = k[X_1,\ldots,X_n,Y_1,\ldots,Y_r]$ $(n \ge 1, r \ge 1)$ be a bigraded polynomial ring with 
$\bideg X_i = (1,0),\ i = 1,\ldots,n,$ and $\bideg Y_j  = (d_i,1),\ j = 1,\ldots,r.$
For $u\ge dv$,  the Hilbert function $H_S(u,v)$ is equal to a polynomial  $$P_S(u,v) = \sum_{i+j = n+r-2}\frac{e_{i,j}}{i!j!}u^iv^j + \text{\rm lower-degree terms}$$
of total degree $n+r-2$ with 
$$e_{i,n+r-2-i} = \left\{\begin{array}{ll} (-1)^{n-i-1}\displaystyle \sum_{j_1+ \ldots + j_r = n-1-i}d_1^{j_1}\cdots d_r^{j_r} & \text{if }\ i < n,\\
0 & \text{if }\ i \ge n.
\end{array}\right.$$
\end{Proposition}

\noindent{\it Remark.} The range $u \ge dv$ can not be removed.
In fact, for $u < d$ and $v=1$, the Hilbert function $H_S(u,1)$ can be equal to different polynomials than $P_S(u,1)$, depending on the range of $u$. For instance,  $H_S(u,1) = 0$ for $u < \min\{d_1,\ldots,d_r\}$.  

Proposition \ref{polynomial} follows from the following result which is based on some polynomial identities of [Ve1]. This result will be used also in Section 5.

\begin{Lemma} \label{observation} Let 
$$f(t) = \frac{e}{m!}t^m + \text{\rm lower-degree terms}$$ 
be a polynomial of degree $m$. Let $d_1,\ldots,d_r$ be a sequence of non-negative integers and $d = \max\{d_1,\ldots,d_r\}$. Set
$$H(u,v) := \sum_{\alpha_1+\ldots+\alpha_r = v} f(u-d_1\alpha_1-\ldots-d_r\alpha_r).$$
Then $H(u,v)$ is equal to a polynomial of degree $m+r-1$. Moreover, if $H(u,v)$ is written in the form
$$H(u,v) = \sum_{i+j=m+r-1} \frac{e_{i,j}}{i!j!}u^iv^j + \text{\rm lower-degree terms},$$
then $$e_{i,m+r-1-i} = \left\{\begin{array}{ll} (-1)^{m-i}e\displaystyle \sum_{j_1+ \ldots + j_r = m-i}d_1^{j_1}\cdots d_r^{j_r} & \text{if }\ i \le m,\\
0 & \text{if }\ i > m.
\end{array}\right.$$
\end{Lemma}

\begin{pf}
For $r=1$ we have 
\begin{align*}
H(u,v) & = \frac{e}{m!}(u-d_1v)^m +\text{\rm lower-degree terms}\\
& = e\sum_{i=0}^m(-1)^{m-i}\frac{d_1^{m-i}}{i!(m-i)!}u^iv^{m-i}+ \text{\rm lower-degree terms}.
\end{align*}
Hence $e_{i,m-i} =  (-1)^{m-i}ed_1^{m-i}$ for $i = 0,\ldots,m$. 
Since $e_{m,0} = e \neq 0$, we get $\deg H(u,v) = m$.  

For  $r>1$ we introduce a new function 
$$H'(u',v') := \sum_{\alpha_1+\ldots +\alpha_{r-1}=v'} f(u'-d_1\alpha_1-\ldots-d_{r-1}\alpha_{r-1}).$$
By induction on $r$ we may assume that 
$$H'(u',v') = \sum_{i+j = m+r-2}\frac{e'_{i,j}}{i!j!}(u')^i(v')^j + \text{\rm lower-degree terms}$$
is a polynomial of degree $m+r-2$. Then
\begin{align*} 
H(u,v) & = \sum_{\alpha_r = 0}^vH'(u-d_r\alpha_r,v-\alpha_r)\\
& = \sum_{\alpha_r = 0}^v\big[\sum_{i+j = m+r-2}\frac{e'_{i,j}}{i!j!}(u-d_r\alpha_r)^i(v-\alpha_r)^j + \text{\rm lower-degree terms}\big].
\end{align*}
Therefore, $H(u,v)$ is a polynomial of degree $m+r-1$ if we can show that
for any pair of non-negative integers $i,j$, the function 
$$G(u,v) = \sum_{\alpha_r=0}^v(u-d_r\alpha_r)^i(v-\alpha_r)^j$$
is a polynomial of degree $i+j+1$.  

We have
\begin{align*}
G(u,v) & = \sum_{\alpha_r=0}^v\sum_{p=0}^i(-1)^p{i \choose p}u^{i-p}(d_r\alpha_r)^p\sum_{q=0}^j(-1)^q{j \choose q}v^{j-q}\alpha_r^q\\
& = \sum_{p=0}^i\sum_{q=0}^j(-1)^{p+q}{i \choose p}{j \choose q}d_r^pu^{i-p}v^{j-q} \sum_{\alpha_r=0}^v\alpha_r^{p+q}.
\end{align*}
By [Ve1, Lemma 2.8] $\sum_{\alpha_r=0}^v\alpha_r^{p+q}$ is a polynomial in $v$ of degree $p+q+1$ with 
$$\sum_{\alpha_r=0}^v\alpha_r^{p+q} = \frac{1}{p+q+1}v^{p+q+1} + \text{\rm lower-degree terms}.$$
Therefore, $G(u,v)$ is a polynomial and we may write
$$ G(u,v) =  \sum_{p=0}^i\sum_{q=0}^j(-1)^{p+q}{i \choose p}{j \choose q}\frac{1}{p+q+1}d_r^pu^{i-p}v^{j+p+1} + \text{\rm lower-degree terms}.$$
By [Ve1, Lemma 2.7] we know that
$$\sum_{q=0}^j(-1)^q{j \choose q}\frac{1}{p+q+1} = \frac{1}{(p+j+1){p+j\choose p}}.$$
This implies
\begin{align*} G(u,v) & =  \sum_{p=0}^i(-1)^{p}{i \choose p}\frac{1}{(p+j+1){p+j\choose p}}d_r^pu^{i-p}v^{j+p+1} + \text{\rm lower-degree terms}\\
& = \sum_{p=0}^i(-1)^{p}\frac{i!j!}{(i-p)!(p+j+1)!}d_r^pu^{i-p}v^{j+p+1} + \text{\rm lower-degree terms}
\end{align*}
The coefficient of $u^iv^{j+1}$ is equal to $1/(j+1)$. Hence $\deg G(u,v) = i+j+1$.\par

Now we are going to compute the coefficients of the highest degree terms of the polynomial $H(u,v)$.
Using the last formula for $G(u,v)$ we have
\begin{align*} 
H(u,v) & = \sum_{i+j = m+r-2}\frac{e'_{i,j}}{i!j!}\sum_{p=0}^i(-1)^{p}\frac{i!j!}{(i-p)!(p+j+1)!}d_r^pu^{i-p}v^{j+p+1}\\
& \hspace{5cm} + \text{\rm lower-degree terms}\\
& = \sum_{i+j = m+r-2}\sum_{p=0}^i(-1)^{p}\frac{e'_{i,j}}{(i-p)!(p+j+1)!}d_r^pu^{i-p}v^{j+p+1} \\
& \hspace{5cm} + \text{\rm lower-degree terms}
\end{align*}
Putting $j = m+r-2-i$ and $h = i-p$ we get
\begin{align*} 
H(u,v) & = \sum_{i=0}^{m+r-2}\sum_{h=0}^i(-1)^{i-h}\frac{e'_{i,m+r-2-i}}{h!(m+r-1-h)!}d_r^{i-h}u^hv^{m+r-1-h} \\
& \hspace{5cm} + \text{\rm lower-degree terms}\\
& = \sum_{h=0}^{m+r-2}\sum_{i=h}^{m+r-2}(-1)^{i-h}\frac{e'_{i,m+r-2-i}}{h!( m+r-1-h)!}d_r^{i-h}u^hv^{m+r-1-h} \\
& \hspace{5cm} + \text{\rm lower-degree terms}
\end{align*}
Thus, if we write
$$H(u,v) = \sum_{i+j = m+r-1}\frac{e_{i,j}}{i!j!}u^iv^j + \text{\rm lower-degree terms},$$
then 
$$e_{i,m+r-1-i} = \left\{\begin{array}{ll} 
 \displaystyle\sum_{h=i}^{m+r-2}(-1)^{h-i}e'_{h,m+r-2-h}d_r^{h-i} & \text{if } i \le m+r-2,\\
0 & \text{if } i = m+r-1.
\end{array}\right.$$

By the induction hypothesis we have
$$e'_{h,m+r-2-h} = \left\{\begin{array}{ll} 
(-1)^{m-h}e\displaystyle\sum_{j_1+\ldots+j_{r-1}=m-h} d_1^{j_1}\ldots d_{r-1}^{j_{r-1}} & \text{if } h \le m,\\
0 & \text{if } h > m.\end{array}\right.$$
It is easy to check that
\begin{align*}
& \sum_{h=i}^{m+r-2}(-1)^{h-i}\big[(-1)^{m-h}e\sum_{j_1+\ldots+j_{r-1}=m-h} d_1^{j_1}\ldots d_{r-1}^{j_{r-1}} d_r^{h-i}\big]\\ &\quad  = (-1)^{m-i}e\sum_{j_1+\ldots+j_r=m-i} d_1^{j_1}\ldots d_r^{j_r}. 
\end{align*}
Therefore,
$$e_{i,m+r-1-i}  =
 \left\{\begin{array}{ll} 
(-1)^{m-i} e\displaystyle\sum_{j_1+\ldots+j_r=m-i} d_1^{j_1}\ldots d_r^{j_r} & \text{if }\ i \le m,\\
0 & \text{if }\ i > m.
\end{array}\right.$$
This completes the proof of Lemma \ref{observation}.
\end{pf}

\noindent{\it Proof of Proposition \ref{polynomial}.}
Let $k[X] = k[X_1,\ldots,X_n]$. Then $k[X]$ is a standard $\NN$-graded algebra with $\deg X_i = 1$.  For all $(u,v) \in \NN^2$ we have
$$S_{(u,v)}= \bigoplus_{\alpha_1+\ldots+\alpha_r=v} k[X]_{u-d_1\alpha_1-\ldots-d_r\alpha_r}Y_1^{\alpha_1}\ldots Y_r^{\alpha_r}.$$
Hence
\begin{align*}
H_S(u,v) & = \sum_{\alpha_1+\ldots+\alpha_r=v} \dim_k k[X]_{u-d_1\alpha_1-\ldots-d_r\alpha_r}\\
& =  \sum_{\alpha_1+\ldots+\alpha_r=v} {u-d_1\alpha_1-\ldots-d_r\alpha_r+n-1 \choose n-1}.
\end{align*}
Let $f(t)$ denote the polynomial ${t+n-1 \choose n-1}$. For $u \ge dv$ we have
$$H_S(u,v) = \sum_{\alpha_1+\ldots+\alpha_r=v} f(u-d_1\alpha_1-\ldots-d_r\alpha_r).$$ 
Hence the conclusion follows from Lemma \ref{observation}.
\qed\medskip

Now we are ready to prove Theorem \ref{exist}.\medskip

\noindent{\it Proof of Theorem \ref{exist}.} First, we represent $R$ as a bigraded quotient of the bigraded polynomial ring $S$ defined in Proposition \ref{polynomial}. 
For any bigraded $S$-module $M$ we set $H_M(u,v) := \dim_kM_{(u,v)}$.
Let $$0 \To F_m \To \cdots \To F_1 \To F_1 \To R \To 0$$
be a bigraded minimal free resolution of $R$ over $S$. Then 
$$H_R(u,v)=\sum_{i=0}^{m}{(-1)}^i H_{F_i}(u,v).$$
Every free module $F_i$ is a direct sum of modules of the form $S(-a,-b)$ for some non-negative integers $a,b$, where $S(-a,-b)_{(u,v)} = S_{(u-a,v-b)}$ for all $(u,v) \in \ZZ^2$. By Proposition \ref{polynomial}, $H_{S(-a,-b)}(u,v) = H_S(u-a,v-b)$ is a polynomial in $u,v$ for $u \ge a+\max\{0,d(v-b)\}$ and $v \ge b$. 
Thus, if we define $u_0 $ to be the maximum of all numbers $a$ and $v_0$ to be the maximum of all numbers $b$, where $(a,b)$ runs all shifting degree occuring in the bigraded minimal free resolution of $R$ over $S$, then the functions $H_{F_i}(u,v)$ and therefore $H_R(u,v)$ are equal to polynomials in $u,v$  for $u \ge  dv+u_0$ and $v \ge v_0$.\qed

\section{Degree of Hilbert polynomial}

Throughout this section,  we assume that the bigraded algebra $R$ is generated by $n$ homogeneous elements 
$x_1,\ldots,x_n$ of bidegree $(1,0)$ and $r$ elements $y_1,\ldots,y_r$ of bidegree $(d_1,1),\ldots,(d_r,1)$, 
where $d_1,...,d_r$ are non-negative integers. 

As we have seen in Theorem \ref{exist},  we can associate with $R$ a Hilbert polynomial $P_R(u,v)$. 
The aim of this section is to compute the total degree $\deg P_R(u,v)$ and the degree $\deg_uP_R(u,v)$   
of $P_R(u,v)$ in the variable $u$.

Let $R_{++}$ denote the ideal of $R$ generated by the elements $x_iy_j$, $i = 1,\ldots,n,\ j = 1,\ldots,r$. 
Let $\Proj R$ denote the set of the bigraded prime ideals $P$ of $R$ with $P \not\supseteq R_{++}$. Let
$$\rdim R := \left\{\begin{array}{ll}
1 & \text{if } \Proj R = \emptyset,\\
\max\{\dim R/P|\ P \in \Proj R\} & \text{if } \Proj R \neq \emptyset.
\end{array}\right.$$
Following [KMV] we call $\rdim R$ the {\it relevant dimension} of the bigraded algebra $R$.

The total degree $\deg P_R(u,v)$ can be expressed in terms of the relevant dimension of $R$ as follows.

\begin{Theorem} \label{degree} $\deg P_R(u,v) = \rdim R-2.$
\end{Theorem}

This theorem covers results of Katz, Mandal and Verma for standard bigraded algebras [KMV, Theorem 2.2]. 
There is another formula for $\deg P_R(u,v)$ given by P. Roberts in [Ro1, Theorem 8.3.4], namely,
$$\deg P_R(u,v) = \max\{\dim R_{(P)}|\ P \in \Proj R\},$$
where  $R_{(P)}$ denotes the degree zero part of the bigraded localization of 
$R$ at $P$ (see [Ro3, Theorem 2(1)] for the case $R$ is generated  
by elements of bidegree $(1,0),(0,1)$, $(1,1)$).
Since the relationship between these two formulas is not trivial, we shall present below 
a direct and short proof for Theorem \ref{degree}. For that we shall need the following notion. 

Let $\Delta = \{(cv,ev)|\ v \in \NN\}$, where $c,e$ are two given positive integers. For any bigraded $R$-module $M$ we define
$$M_\Delta := \oplus_{v \in \ZZ}M_{(cv,ev)}.$$
It is clear that $R_\Delta$ is a $\NN$-graded algebra and $M$ a graded $R_\Delta$-module. We call $R_\Delta$  a {\it diagonal subalgebra} of $R$. Diagonal subalgebras were introduced in [STV], [CHTV]. They have been studied mainly in the case $R$ is the Rees algebra of a homogeneous ideal.

\begin{Lemma} \label{diagonal} 
Let $d= \max\{d_1,\ldots,d_r\}$. Then\par
{\rm (a) }  $R_\Delta$ is a standard graded algebra if $c \ge de$,\par
{\rm (b) }  $\dim R_\Delta = \rdim R - 1$ if $c > de$.
\end{Lemma}

\begin{pf}   
(a) By definition, $(R_\Delta)_v = R_{(cv,ev)}$. Therefore, $(R_\Delta)_v$ is the vector space spanned by the products $x_1^{a_1}\cdots x_n^{a_n}y_1^{b_1}\cdots y_r^{b_n}$ of bidegree $(cv,ev)$, which means 
\begin{align*}
a_1 + \ldots + a_n + b_1d_1 + \ldots + b_rd_r & = cv,\\
b_1 + \ldots + b_r & = ev.
\end{align*}
If $v \ge 2$, a product $x_1^{a_1}\cdots x_n^{a_n}y_1^{b_1}\cdots y_r^{b_n}$ of bidegree $(cv,ev)$ is always divisible by a product $x_1^{a'_1}\cdots x_n^{a'_n}y_1^{b'_1}\cdots y_r^{b'_n}$ of bidegree $(c(v-1),e(v-1))$.
Indeed, we first choose non-negative integers $b'_1 \le b_1,\ldots,b'_r \le b_r$  such that 
$$b'_1 + \ldots + b'_r = e(v-1).$$
Since $(b_1-b'_1) + \ldots + (b_r-b'_r) = e$, we have 
$$c \ge de \ge (b_1-b'_1)d_1 + \ldots + (b_r-b'_r)d_r .$$
From this it follows that
$$c(v-1) - (b'_1d_1 + \ldots + b'_rd_r) \le cv - (b_1d_1 + \ldots + b_rd_r) = a_1+\ldots+a_n.$$
Hence there exist non-negative integers $a'_1 \le a_1,\ldots,a'_n \le a_r$ such that
$$a'_1+\ldots+a'_n = c(v-1)- b'_1d_1 + \ldots + b'_rd_r.$$
Clearly, the chosen integers $a'_1,\ldots,a'_n,b'_1,\ldots,b'_r$ yield a product $x_1^{a'_1}\cdots x_n^{a'_n}y_1^{b'_1}\cdots y_r^{b'_n}$ as required. So we can conclude that the $\NN$-graded algebra $R_\Delta$ is generated by elements of degree 1.\par

(b) It is easy to check that if $P$ is a bigraded prime ideal of $R$ and $I$ is a bigraded $P$-primary ideal, then $P_\Delta$ is a prime ideal of $R_\Delta$ and $I_\Delta$ is a graded $P_\Delta$-primary ideal. 
Further, if $0_R = \cap I_j$ is a graded primary decomposition of the zero ideal $0_R$ of $R$, then $0_{R_\Delta} = \cap(I_j)_\Delta$ is a graded primary decomposition of the zero ideal $0_{R_\Delta}$ of $R_\Delta$. Thus, to show that $\dim R_\Delta = \rdim R-1$ we may assume that
$R$ is a domain. \par

If $\Proj R = \emptyset$, we have $\rdim R = 1$ and $(R_{++})^v = 0$ for some positive integer $v$. 
Since from the assumption $c > de$ we can deduce that $R_{(cv,ev)} \subseteq (R_{++})^v$, we get $(R_\Delta)_v = 0$. Therefore, $\dim R_\Delta = 0$. 

If $\Proj R \neq \emptyset$, the zero ideal belongs to $\Proj R$. This implies $R_{++} \neq 0$ and $\rdim R = \dim R$. In this case, we will use the formulas
\begin{align*}
\dim R & = \td_kQ(R),\\
\dim R_\Delta &= \td_kQ(R_\Delta),
\end{align*}
where $Q(R)$ and $Q(R_\Delta)$ denote the quotient fields of $R_\Delta$ and $R$, respectively. Since $R_{++} \neq 0$, we can find  a homogeneous element $x \in R$ of bidegree $(1,0)$. Since $e \ge 1$, $x \not\in R_\Delta$. Because of its degree, the element $x$ can not be a root of any algebraic equation over $Q(R_\Delta)$. Hence $x$ is a transcendent element over $Q(R_\Delta)$. It suffices now to show that $Q(R)$ is an algebraic extension of $Q(R_\Delta)(x)$. For this will imply
$$\td_kQ(R) = \td_kQ(R_\Delta) +1.$$

Since $R_{++} \neq 0$, there exists an element $y\in R$ of bidegree $(d,1)$. For any generator $x_i$ of $R$ we can write $x_i = x_ix^{c-de}y^e/x^{c-de}y^e$.   Since $\bideg x_ix^{c-1-d_1e}y^e = \bideg x^{c-de}y^e = (c,e)$, we get $x_i \in Q(R_\Delta)(x)$. 
For any generator $y_j$ of $R$ we have $\bideg x^{c-d_je}y_j^e = (c,e).$
As $c > d_je$, we get $x^{c-d_je}y_j^e \in R_\Delta$. Hence $y_j$ is an algebraic element over $Q(R_\Delta)(x)$. 
Since $R$ is generated by the elements $x_i$ and $y_j$, we can conclude that $Q(R)$ is an algebraic extension of $Q(R_\Delta)(x)$. 
\end{pf}

\noindent{\it Proof of Theorem \ref{degree}.}
Let $s$ be the degree of the Hilbert polynomial $P_R(u,v)$ and write 
$$P_R(u,v)=\sum_{i+j=s}\frac {e_{i,j}}{i!j!}u^iv^j+\text{ lower-degree terms}.$$
We consider the diagonal subalgebras $R_\Delta$ for any $c > d$ and $e=1$. Let 
$$H_{R_\Delta}(v) := H_R(cv,v) = \dim_k R_{(cv,v)}$$
for all $v \in \NN$.
Then $H_{R_\Delta}(v)$ is the Hilbert function of $R_\Delta$.
By Lemma \ref{diagonal},  $H_{R_\Delta}(v)$ is equal to a polynomial $P_{R_\Delta}(v)$ of degree $\rdim R-2$ for $v$ large enough. Since $cv \ge dv+v$, by Theorem \ref{exist} we also have $H_R(cv,v) = P_R(cv,v)$ for $v$ large enough.   Therefore,
$$P_{R_\Delta}(v) = P_R(cv,v) = \sum_{i=0}^s\frac{e_{i,s-i}}{i!(s-i)!}c^iv^s+\text{ lower-degree terms}.$$
If we choose $c$ large enough then 
$\sum_{i=0}^s\frac{e_{i,s-i}}{i!(s-i)!}c^i \neq 0.$ 
So we can conclude that $s = \rdim R-2$. \qed\medskip

Now we shall see that $\deg P_R(u,v)$ can be expressed as the dimension of a quotient ring of $R$.  \par

Let $\bar R := R/0:(R_{++})^\infty,$ where $0:(R_{++})^\infty$ denotes the  ideal of the elements of $R$ annihilated by some power of $R_{++}$. Note that $0:(R_{++})^\infty$ can be obtained from any decomposition of the zero ideal of $R$ by deleting those primary components whose associated primes contain $R_{++}$. 
Then $\rdim R = \dim \bar R$ if $\Proj R \neq \emptyset$, and $\bar R = 0$ if $\Proj R = \emptyset$. Thus, we always have
$$\deg P_R(u,v) = \dim \bar R -2.$$

The degree $\deg_u P_R(u,v)$ can be also expressed in terms of $\bar R$ by relating $P_R(u,v)$ to the Hilbert polynomials of the graded modules
$$R_v := \oplus_{u \ge 0}R_{(u,v)}.$$
Note that $R_0$ is a finitely generated standard $\NN$-graded algebra
and $R_v$ is a finitely generated graded $R_0$-module. We set 
$$\sdim R := \dim \bar R_0$$
and call it the $x$-{\it dimension} of $R$. The name stems from [Ro3] where $\sdim R$ is defined to be 
the dimension of $R_v$ for $v$ large enough. But the following formulas for  $\deg_u P_R(u,v)$ show 
that both definitions are equal. These formulas were already proved for standard bigraded algebras 
in [Tr2, Theorem 1.7] and implicitly for bigraded algebras generated by elements of bidegree $(1,0),(0,1),(1,1)$ 
in [Ro3, Theorem 2(2)]. 

\begin{Theorem} \label{partial} 
$\deg_u P_R(u,v)  = \dim R_v-1 = \sdim R-1$ for $v$ large enough.
\end{Theorem} 

\begin{pf}  Let $H_{R_v}(u)$  and $P_{R_v}(u)$ denote the Hilbert function and the Hilbert polynomial of $R_v$. 
Then $H_{R_v}(u) = \dim_kR_{(u,v)}$ for all $u \in \NN$. By Theorem \ref{exist}, there exist integers $u_0,v_0$ 
such that $H_{R_v}(u) = P_R(u,v)$ for $u \ge cv+u_0$ and $v \ge v_0$. Therefore $P_{R_v}(u) = P_R(u,v)$ for $v \ge v_0$. For $v$ large enough, the leading coefficient of $P_R(u,v)$ as a polynomial in $u$ does not vanish. From this we deduce that $\deg_uP_R(u,v) = \deg P_{R_v}(u)$.
Since $\deg P_{R_v}(u) = \dim R_v-1,$ that implies 
$$\deg_uP_R(u,v) = \dim R_v-1.$$

Choose a positive integer $m$ such that $0:(R_{++})^\infty = 0:(R_{++})^m$. For $u \ge cv+m$  and $v \ge m$,  
we have $R_{(u,v)} \subset (R_{++})^m$; hence
$$\big(0:(R_{++})^\infty\big)_{(u,v)} \subseteq  (R_{++})^m \big(0:(R_{++})^\infty\big) = 0,$$
which implies $H_R(u,v) = H_{\bar R}(u,v)$. Therefore,
$P_R(u,v) = P_{\bar R}(u,v)$. As we have seen above, $\deg_u P_{\bar R}(u,v) = \dim \bar R_v-1$. But $\dim \bar R \le \sdim R$. So we get 
$$\deg_u P_R(u,v) = \deg_uP_{\bar R}(u,v) \le \sdim R-1.$$
 
It remains to show that $\dim R_v \ge  \sdim R$ for $v$ large enough.
For this will imply the converse inequality $\deg_uP_R(u,v) \ge \sdim R-1$; hence $\deg_uP_R(u,v) = \sdim R-1$. 
Let $R_+$ denote the ideal $\oplus_{v > 0}R_v$ of $R$. Then $R_v \cong (R_+)^v/(R_+)^{v+1}$
and $\bar R_0 = R/(R_++0:(R_{++})^\infty)$. Let $P$ be an associated prime of $R_++0:(R_{++})^\infty$ such that $\dim R/P = \dim \bar R_0 = \sdim R$.   Then $0:(R_+)^v \subseteq 0: (R_{++})^\infty \subset P.$
Therefore, $(0:(R_+)^v)_P$ is a proper ideal in the local ring $R_P$ so that $(R_+)^v_P \neq 0$. By Nakayama's lemma, this implies 
$(R_+)^v_P/(R_+)^{v+1}_P \neq 0$.  Hence
$$\dim R_v = \dim (R_+)^v/(R_+)^{v+1} \ge \dim R/P = \dim \bar R_0.$$
\end{pf}

It is easy to find examples with $\dim  R_0 > \dim \bar R_0 = \sdim R$.

\begin{Example} 
{\rm Let $R = k[X_1,X_2,Y_1]/(X_1Y_1)$ with $\bideg X_1 = \bideg X_2 = (1,0)$ and $\bideg Y_1 = (1,1)$. Then $\bar R = k[X_1,X_2,Y_1]/(X_1) = k[X_2,Y_1]$. Since $R_0 = k[X_1,X_2]$ and $\bar R_0 = k[X_2]$, we have
$\dim R_0 = 2 > 1 = \dim \bar R_0$.} \end{Example}

By Theorem \ref{degree} and Theorem \ref{partial} we always have
$$\rdim R = \deg P_R(u,v) \ge \deg_u P_R(u,v) = \sdim R + 1.$$
The following example shows that the inequality may be strict.

\begin{Example}
{\rm Let $R = k[X_1,X_2,Y_1,Y_2,Y_3]/(X_1Y_1)$ with $\bideg X_i =(1,0)$ and $\bideg Y_j = (d_j,1)$ for any sequence of non-negative integers $d_1,d_2,d_3$. Then $\bar R = R$ and $R_0 = k[X_1,X_2]$. Hence
$\rdim R = 4 > 2+1 = \sdim R+1$.}
\end{Example}

We are not able to find a formula for $\deg_vP_R(u,v)$ (the degree of $P_R(u,v)$ in $v$). The method of Theorem \ref{partial} does not work in this case since by the range of the equality $H_R(u,v) = P_R(u,v)$ we can not fix $u$ to estimate $\deg_vP_R(u,v)$. 

\section{Leading coefficients of Hilbert polynomial}

Let $R$ be a bigraded algebra $R$ generated by $n$ homogeneous elements $x_1,\ldots,x_n$ of bidegree $(1,0)$ and $r$ elements $y_1,\ldots,y_r$ of bidegree $(d_1,1),\ldots,(d_r,1)$, where $d_1,...,d_r$ are non-negative integers. \par

We write the Hilbert polynomial $P_R(u,v)$ in the form
$$P_R(u,v) = \sum_{i= 0}^s \frac{e_i(R)}{i!(s-i)!}u^iv^{s-i} + \text{\rm lower-degree terms},$$
where $s = \deg P_R(u,v)$. Note that $P_R(u,v) = 0$ if $s=-1$. Following [Te] we call the numbers $e_i(R)$ the {\it mixed multiplicities} of $R$. Moreover, we set
$$\rho_R := \max\{i|\ e_i(R) \neq 0\}.$$
The aim of this section is to establish some basic properties of mixed multiplicities.  

We shall need the following technical notion. Let $x$ be a homogeneous element of $R$. We call $x$ a {\it filter-regular} element of $R$ if $x \not\in P$ for any associated prime ideal $P \in \Proj R$ of $R$. This notion has its origin in the theory of Buchsbaum rings (see e.g. [SV]). It has been already used to study Hilbert polynomials of standard bigraded algebras in [Tr2]. It is obvious that if the base field $k$ is infinite and $R_{(1,0)} \neq 0$ then we can always find a filter-regular element of degree $(1,0)$.\par

\begin{Lemma} \label{filter}
Let $x$ be a filter-regular element of degree $(1,0)$ in $R$. Then\par
{\rm (a) } $\deg P_{R/xR}(u,v) \le \deg P_R(u,v)-1$,\par
{\rm (b) } $\deg_u P_{R/xR}(u,v) = \deg_u P_R(u,v)-1$,\par
{\rm (c) } $e_{i-1}(R/xR) = e_i(R)$ for $i \ge 1$, \par
{\rm (d) } $\rho_{R/xR} = \rho_R-1$ if $\deg P_{R/xR}(u,v) = \deg P_R(u,v)-1$. 
\end{Lemma}

\begin{pf}
Since $x \not\in P$ for any associated prime ideal $P \not\supseteq R_{++}$ of $R$,  $0:x \subseteq 0:(R_{++})^\infty$. As we have seen in the proof of Theorem \ref{partial}, there is an integer $m$ such that $(0:x)_{(u,v)} \subseteq \big(0:(R_{++})^\infty\big)_{(u,v)} = 0$ for $u > dv+m$ and $v \ge m$. By Theorem \ref{exist}, this implies $P_R(u,v) = P_{R/0:x}(u,v)$. Rewrite $P_R(u,v)$ in the form
$$P_R(u,v) = \sum_{i+j \le s} e_{i,j}{u \choose i}{v \choose j},$$
and consider the exact sequence of bigraded algebras
$$0 \To R/0:x \overset{x} \To R \To R/xR \To 0.$$
Since $x$ is a homogeneous element of degree $(1,0)$, 
$$\displaylines{P_{R/xR}(u,v) = P_R(u,v) - P_{R/0:x}(u-1,v) = P_R(u,v) - P_R(u-1,v) \cr
= \sum_{i+j \le s}e_{i,j} {u \choose i}{v \choose j}\ - \sum_{i+j \le s}e_{i,j}{u-1 \choose i}{v \choose j}\cr
= \sum_{\tiny \begin{array}{c}i+j \le s \\ i \ge 1\end{array}} \left[{e_{i,j}\over (i-1)!j!}u^{i-1}v^j + \text{terms of degree $< i+j-1$ with $\deg_u < i-1$}\right].}$$
Now we can easily verify (a), (b) and (c). Moreover, we have $\rho_{R/xR} = \rho_R-1$ if $\rho_R > 0$ and $\deg P_{R/xR}(u,v) < \deg P_R(u,v)-1$ if $\rho_R = 0$. Hence (d) is obvious. 
\end{pf}

\begin{Theorem} \label{coefficient}
The mixed multiplicities $e_i(R)$ are integers with $e_{\rho_R}(R) > 0$.
\end{Theorem}

\noindent{\it Remark.} As shown in Proposition \ref{polynomial}, a mixed multiplicity can be a negative number. 
However, if $R$ is a standard bigraded algebra, all mixed multiplicities are non-negative [W]. 
Moreover, $e_0(R),\ldots,e_s(R)$ can be any sequence of non-negative integers with at least one positive entry 
[KMV, Example 5.2]. 

\begin{pf}
We will first show that $e_{\rho_R}(R) > 0$ is an integer. Consider the diagonal subalgebra $R_\Delta := \oplus_{v \ge 0}R_{(cv,v)}$ for any integer $c > \max\{d_1,\ldots,d_r\}$. By Lemma \ref{diagonal}(a),  $R_\Delta$ is a standard graded algebra. By the proof of Theorem \ref{degree} we have
$$P_{R_\Delta}(v) =  \sum_{i=0}^s\frac{e_i(R)}{i!(s-i)!}c^iv^s+\text{ lower-degree terms}$$
with $\deg P_{R_\Delta}(v) = s$. From this it follows  that
$$\sum_{i=0}^s\frac{e_i(R)}{i!(s-i)!}c^i  = \frac{e(R_\Delta)}{s!} > 0.$$
The leading coefficient of the polynomial in $c$ on the left side must be positive. Hence $e_{\rho_R}(R)  > 0$.\par
From the above formula we also get
$$e_0(R) = e(R_\Delta) - \sum_{i=1}^s{s \choose i}e_i(R).$$
Hence $e_0(R)$ is an integer. Since $\deg_u P_R(u,v) \ge \rho_R$,
to show that $e_1(R),\ldots,e_{\rho_R}(R)$ are integers we may assume that $\deg_u P_R(u,v) > 0$. Extending the base field $k$ by an indeterminate, we may also assume that $k$ is an infinite field. Then we can find a filter-regular element $x \in R$ of degree $(1,0)$. 
By Lemma \ref{filter}(b), $\deg_u P_{R/xR}(u,v) = \deg_u P_R(u,v)-1$. Using induction on $\deg_u P_R(u,v)$ we may assume that $e_i(R/xR)$ is an integer for $i \ge 0$. By Lemma \ref{filter}(c), this implies that $e_i(R)$ is an integer for $i \ge 1$. 
\end{pf}

Now we will give a condition for $\rho_R = \deg_uP_R(u,v)$ (see Example \ref{less} for an example with 
$\rho_R < \deg_uP_R(u,v)$). Note that this condition is satisfied if $R$ is a domain or a Cohen-Macaulay ring.

\begin{Proposition} \label{equi}
Assume that $\dim R/P = \rdim R$ for all minimal prime ideals of $\Proj R$.
Then $e_i(R) > 0$ for $i = \deg_u P_R(u,v)$.
\end{Proposition}

\begin{pf} By Theorem \ref{coefficient} we only need to show that $\rho_R = \deg_u P_R(u,v)$. For this we may assume that $\deg_u P_R(u,v) > 0$. By Theorem \ref{partial} we have $R_{(1,0)} \neq 0$.  Extending the base field $k$ by an indeterminate we may assume that $k$ is an infinite field. Then we can find a filter-regular element $x \in R$  of degree $(1,0)$. By Lemma \ref{filter}(b),
$$\deg P_{R/xR}(u,v) \ge \deg_u P_{R/xR}(u,v) = \deg_uP_R(u,v)-1 \ge 0.$$
By Theorem \ref{degree}, this implies $\rdim R/xR \ge 2$. Hence $\Proj R/xR \neq \emptyset$. Let $P$ be an arbitrary minimal associated prime ideal of $xR$ such that $P\not\supseteq xR + R_{++}$. Then $P \in \Proj R$. Let $P' \subseteq P$ be a minimal prime ideal of $\Proj R$. Then $P'$ is an associated prime ideal of $R$. Hence $x \not\in P'$. From this it follows that
$\dim (R/P')_P = 1$.  Using the assumption we obtain
$$\dim R/P = \dim R/P' - \text{ht} P/P' = \rdim R - \dim (R/P')_P = \rdim R-1.$$
So we have proved that $\rdim R/xR = \rdim R-1$ and that $R/xR$ satisfies the assumption of the theorem. By Theorem \ref{degree} we have 
$$\deg P_{R/xR}(u,v) = \rdim R/xR-2 = \rdim R -3 = \deg P_R(u,v)-1.$$
By Lemma \ref{filter}(d), this implies $\rho_{R/xR} = \rho_R-1$. By induction we may assume that 
$$\rho_{R/xR} = \deg_uP_{R/xR}(u,v) = \deg_uP_R(u,v)-1.$$
Hence we can conclude that $\rho_R = \deg_uP_R(u,v)$.
\end{pf} 

Like for the usual multiplicity, we have an associativity formula for mixed multiplicities of bigraded algebras.

\begin{Proposition} \label{associative}
Let ${\cal A}(R)$ be the set of the associated prime ideals $P \in \Proj R$ with $\dim R/P = \rdim R$. For any $P \in {\cal A}(R)$ let $\ell(R_P)$ denote the length of the local ring $R_P$. Then
$$e_i(R)= \sum_{P \in {\cal A}(R)}\ell(R_P)e_i(R/P).$$
\end{Proposition}

\begin{pf} By [Ma, Theorem 6.4] there exists a filtration
$$R = Q_1 \supset Q_2 \supset ... \supset Q_m \supset Q_{m+1} = 0$$
of bigraded ideals of $R$ such that $Q_j/Q_{j+1} \cong (R/P_j)(a_j,b_j)$ for some associated prime ideal $P_j$ of $R$ and integers $a_j,b_j$, $j = 1,...,m$.  From this it follows that
$H_R(u,v) = \sum_{j=1}^mH_{R/P_j}(u+a_j,v+b_j).$ Hence
$P_R(u,v) = \sum_{j=1}^mP_{R/P_j}(u+a_j,v+b_j).$ By Theorem \ref{degree} this implies
$$e_i(R) = \sum_{P_j \in {\cal A}(R)}e_i(R/P_j).$$
For every prime ideal $P \in {\cal A}(R)$ we consider the induced filtration
$$R_P = (Q_1)_P \supseteq (Q_2)_P \supseteq ... \supseteq (Q_m)_P \supseteq (Q_{m+1})_P = 0.$$
We have $(Q_j/Q_{j+1})_P = 0$ if $P_j \neq P$ and $(Q_j/Q_{j+1})_P \cong k(P)$ if $P_j = P$, where $k(P)$ denotes the residue field of $R_P$. Therefore, $\ell(R_P)$ is the number of indices $j$ for which $P_j = P$. So we can conclude that
$$e_i(R) = \sum_{P \in {\cal A}(R)}\ell(R_P)e_i(R/P).$$
\end{pf}

By Proposition \ref{associative}, the mixed multiplicities of $R$ depend only on the prime ideals of $\Proj R$ with the highest dimension. Using this fact we can easily construct examples with $\rho_R < \deg_uP_R(u,v)$. 

\begin{Example} \label{less}
{\rm Let $R = k[X_1,X_2,Y_1,Y_2,Y_3]/(X_1)\cap(Y_1,Y_2)$ with $\bideg X_i = (1,0)$ and $\bideg Y_j = (d_j,1)$, where $d_1,d_2,d_3$ can be any sequence of non-negative integers. Then 
$$\deg_u P_R(u,v) = \sdim R-1 = \dim k[X_1,X_2] - 1 = 1.$$
By Proposition \ref{associative}, the mixed multiplicities of $R$ are equal to those of the quotient ring $R/X_1R$. Hence 
$$\rho_R = \rho_{R/X_1R} \le \deg_u P_{R/X_1R}(u,v) =\sdim R/X_1R -1 = \dim k[X_2]-1 = 0.$$}
\end{Example}

Combining Proposition \ref{equi} with Proposition \ref{associative} we obtain the following generalization of a result of P. Roberts in the case
$R$ is generated by elements of degree $(1,0),(0,1),(1,1)$  [Ro3, Theorem 2(3)].
This result was used to give a criterion for the positivity of Serre's intersection multiplicity [Ro3, Proposition 6].

\begin{Corollary} Assume that there exists a prime ideal $P \in {\cal A}(R)$ with $\sdim R/P = \sdim R$. Then $e_i(R) > 0$ for $i = \sdim R-1$.
\end{Corollary}

\begin{pf}  
Let $m := \sdim R$. Then $\sdim R/P \le m$ for every prime ideal $P \in {\cal A}(R)$. 
If $\sdim R/P < m$, then
$\deg_uP_{R/P}(u,v) < m-1$. Hence $e_{m-1}(R/P) = 0$. If $\sdim R/P = m$, then $\deg_uP_{R/P}(u,v) = m-1$. Hence $e_{m-1}(R/P) > 0$ by Proposition \ref{equi}. Applying Proposition \ref{associative} we get 
$$e_{m-1}(R) =\sum_{P \in {\cal A}(R)} \ell(R_P)e_{m-1}(R/P) > 0.$$
\end{pf}

\section{Rees algebras of homogeneous ideals}

Let $A$ be a standard graded algebra over a field $k$. Let $I$ be a homogeneous ideal of $A$. The {\it Rees algebra} $A[It]$ of $I$ is the subring $\oplus_{v \ge 0}I^vt^v$ of $A[t]$. It has a natural bigraded structure by setting
$$A[It]_{(u,v)} := (I^v)_ut^v$$
for all $(u,v) \in \NN^2$. 

Let $\mm = (x_1,\ldots,x_n)$ be the maximal graded ideal of $A$, where $x_1,\ldots,x_n$ are homogeneous elements with $\deg x_i = 1$. Let $I = (f_1,\ldots,f_r)$, where $f_1,\ldots,f_r$ are homogeneous elements with $\deg f_j = d_j$. Put $y_j = f_jt$. Then $A[It]$ is generated by the elements $x_1,\ldots,x_n$ and $y_1,\ldots,y_r$ with $\bideg x_i = (1,0)$ and $\bideg y_j = (d_j,1)$. Hence $A[It]$ belongs to the class of bigraded algebras considered in the preceding section. \par

Let $0:I^\infty$ denote the ideal of the elements of $A$ 
which are annihilated by some power of $I$.  Note that $\dim A/0:I^\infty = \dim A$ if $I$ has positive height.

\begin{Lemma} \label{dim-Rees}
Let $A[It]$ be the above bigraded Rees algebra. Then\par
{\rm (a) } $\rdim A[It] = \dim A/0:I^\infty+1$,\par
{\rm (b) } $\sdim A[It]  = A/0:I^\infty$.
\end{Lemma}

\begin{pf} (a) For every primary ideal $Q$ of $A$ let $Q^* := \oplus_{v\ge 0}(Q\cap I^v)t^v$. It is easy to see that if $P$ is a prime ideal of $A$ and $Q$ a $P$-primary ideal, then $P^*$ is a prime ideal of $A[It]$ and $Q^*$ a $P^*$-primary ideal. Let $A[It]_{++}$ be the ideal of $A[It]$ generated by the elements $x_iy_j$, $i = 1,\ldots,n$, $j = 1,\ldots,r$. Then $A[It]_{++} = \oplus_{v \ge 1}\mm I^vt^v$. Therefore, $P^* \not\supseteq A[It]_{++}$ if and only if $P \not\supseteq I$. By [Va, Section 1] we have
$$\dim A[It]/P^* =  \dim A/P+1$$
if $P \not\supseteq I$. Let $0_A = \cap Q_i$ be a primary decomposition of the zero ideal of $A$.
Then $0_{A[It]} = \cap Q_i^*$ is a primary decomposition of the zero ideal of $A[It]$.  Therefore, $P^*$ is an associated prime of $A[It]$ with $P^* \not\supseteq A[It]_{++}$ if and only if $P$ is an associated prime ideal of $A$ with $P\not\supseteq I$, which is equivalent to the condition that $P$ is an associated prime of $A/0:I^\infty$. 

If $\Proj A[It] = \emptyset$, then $I \subseteq P$ for every associated prime ideal $P$ of $A$. Hence $I \subseteq \sqrt{0_A}$ so that $0:I^\infty = A$. In this case, we have
$$\rdim A[It] = 1 = 0+1 =\dim A/0:I^\infty+1.$$ 
If $\Proj A[It] \neq \emptyset$, the above analysis implies 
 \begin{align*}
\rdim A[It] & = \max\{\dim A[It]/P^*|\ P \in \Ass A,\ P\not\supseteq I\} \\
& = \max\{\dim A/P+1|\ P \in \Ass A, P\not\supseteq I\}\\
& = \dim A/0:I^\infty +1.
\end{align*} \par

(b) It is easy to see that $A[It]_0 = A$ and $A[It]_v \cong I^v$ for $v \ge 1$. 
For $v$ large enough, we have $\text{ann}_A(I^v) = 0:I^v = 0:I^\infty$, hence
$$\dim A[It]_v = \dim I^v = \dim A/0:I^\infty.$$
By Theorem \ref{partial}, this implies $\sdim A[It] = \dim A/0:I^v = \dim A/0:I^\infty.$
\end{pf}

For any finitely generated graded $A$-module $M$ we will denote by $e(M)$, $H_M(u)$ and $P_M(u)$ the multiplicity, the Hilbert function and the Hilbert polynomial of $M$.

We can summarize the results of the preceding sections for the bigraded Rees algebra $A[It]$ as follows.

\begin{Theorem} \label{Rees} 
Set $d = \max\{d_1,\ldots,d_r\}$
and $s = \dim A/0:I^\infty-1$.  There exist  integers $u_0,v_0$ such that for  
$u \ge dv+u_0$ and $v \ge v_0$, the Hilbert function $H_{A[It]}(u,v)$ is equal to a polynomial $P_{A[It]}(u,v)$ with 
$$\deg P_{A[It]}(u,v) = \deg_u P_{A[It]}(u,v) = s.$$
Moreover, if $s \ge 0$ and $P_{A[It]}(u,v)$ is written in the form 
$$P_{A[It]}(u,v) = \sum_{i=0}^s\frac{e_i(A[It])}{i!(s-i)!}u^iv^{s-i}+  \text{\rm lower-degree terms},$$
then the coefficients $e_i(A[It])$ are integers for all $i$ with $e_s(A[It])  = e(A/0:I^\infty)$. \end{Theorem}

\begin{pf} Except the last formula for $e_s(A[It])$, all statements of Theorem \ref{Rees} follow from 
Theorem \ref{exist}, Theorem \ref{degree}, Theorem \ref{partial}, and Theorem \ref{coefficient} by taking Lemma \ref{dim-Rees} into consideration. \par

To prove the formula for $e_s(A[It])$ we first note that $H_{A[It]}(u,v)=H_{I^v}(u)$ for all $u$. By Theorem \ref{exist}, this implies 
$$ P_{A[It]}(u,v)=P_{I^v}(u) $$ 
for $v$ large enough. Since
$\deg P_{A[It]}(u,v) = \deg_u P_{A[It]}(u,v) = s$, the term of $u^s$ in 
the polynomial $P_{A[It]}(u,v)$ must be non-zero. The coefficient of this term is
equal to $\frac{e_s(A[It])}{s!}$. Therefore, for $v$ large enough, we may write
$$P_{I^v}(u) = \frac{e_s(A[It])}{s!}u^s +  \text{\rm lower-degree terms}.$$
Since $\dim I^v = s+1$, this implies $e_s(A[It]) = e(I^v)$. Note that $I^v \cap (0:I^\infty) = 0$.
Then $I^v \cong I^v+(0:I^\infty)/0:I^\infty$. Hence there is an exact sequence of the form
$$0 \To I^v  \To A/0:I^\infty \To A/I^v+(0:I^\infty) \To 0.$$
Since $(0:I^\infty):I = 0:I^\infty$, we have 
$\dim A/0:I^\infty > \dim A/I^v+(0:I^\infty)$. Therefore, the above sequence implies 
$e(I^v) = e(A/0:I^\infty)$. So we can conclude that $e_s(A[It]) =  e(A/0:I^\infty)$.
\end{pf}

\noindent{\it Remark.} There is another kind of mixed multiplicities associated with $I$. Let $M:=(\mm,It)$. The associated graded ring $R := \oplus_{n \ge 0}{M^n}/{M^{n+1}}$ has a natural standard bigrading with 
$R_{(u,v)} := {\mm^uI^v}/{\mm^{u+1}I^v}.$ Teissier [Te] called  $e_i(R)$ a mixed multiplicity of the pair $(\mm,I)$ and denoted it by $e_i(\mm|I)$. The multiplicity of the Rees algebra $A[It]$ and of the extended Rees algebra $A[It,t^{-1}]$ can be expressed in terms of the mixed multiplicities $e_i(\mm|I)$ (see [Ve1], [Ve2], [KV]). One may use $e_i(A[It])$ to compute $e_i(\mm|I)$ if $I$ is generated by forms of the same degree. In this case, there is a bigraded isomorphism $R \cong A[It]$ with a linear transformation of the bidegree.  

The Hilbert polynomial $P_{A[It]}(u,v)$ can be used to compute
the Hilbert polynomial of the quotient ring $A/I^v$ for $v$ large enough.

\begin{Corollary} \label{power} For $v$ large enough,
$$P_{A/I^v}(u) = P_A(u) - P_{A[It]}(u,v).$$
\end{Corollary}

\begin{pf} 
By the proof of  Theorem \ref{Rees} we have $P_{I^v}(u) = P_{A[It]}(u,v)$ for $v$ large enough. Hence the conclusion is immediate.
\end{pf}

Let  $V$ denote the blow-up of the subscheme of $\Proj A$ defined by $I$. 
It is known that $V$ can be embedded into a projective space by the linear system $(I^e)_c$ for any pair of positive integers $e,c$ with $c > de$ [CH, Lemma 1.1], and that the embeddings often give concrete varieties with interesting algebraic properties (see e.g.~[Gi], [GG], [GGH], [GL]). Let $V_{c,e}$ denote the embedded variety. We can use the Hilbert polynomial $P_{A[It]}(u,v)$ to describe the degree of $V_{c,e}$ as follows.

\begin{Corollary} \label{embedded} Let $s = \dim A/0:I^\infty-1$. Assume that $c > de$. Then
$$\deg V_{c,e} = \displaystyle \sum_{i=0}^s {s \choose i}e_i(A[It])c^ie^{s-i}.$$ 
\end{Corollary}

\begin{pf}   The homogeneous coordinate ring of the embedded variety $V_{c,e}$ is the subalgebra $k[(I^e)_c]$ of $A$ generated by the elements of $(I^e)_c$.  As observed in [STV] and [CHVT], we may identify $k[(I^e)_c]$ with the diagonal subalgebra $A[It]_\Delta$ of the bigraded Rees algebra $A[It]$, where 
$\Delta =\{(cv,ev)|\ v \in \NN\}$. Hence 
$$\deg V_{c,e}  = e(A[It]_\Delta).$$
Let $P_{A[It]_\Delta}(v)$ denote the Hilbert polynomial of $A[It]_\Delta$. As we have seen in the proof of Theorem \ref{degree}, 
$$P_{A[It]_\Delta}(v) = P_{A[It]}(cv,v).$$
By Theorem \ref{Rees} we may write 
$$ P_{A[It]}(cv,v) = \sum_{i=0}^s\frac{e_i(A[It])}{i!(s-i)!}c^ie^{s-i}v^s +
\text{\rm lower-degree terms}.$$
By Lemma \ref{diagonal}(b) and Lemma \ref{dim-Rees}(a),
$$\deg P_{A[It]_\Delta}(u) = \dim A[It]_\Delta-1 = \rdim A[It]-2 = s.$$
Hence we can conclude that
$$e(A[It]_\Delta) = \sum_{i=0}^s {s \choose i} e_i(A[It]) c^ie^{s-i}.$$
\end{pf}

Now we will present a method for the computation of $H_{A[It]}(u,v)$. This method  was introduced in [HeTU] in order to compute the multiplicity of $A[It]$ (see also [RaS], [Tr1]).

Let $S=A[Y_1,\ldots,Y_r]$ be a polynomial ring over $A$. Mapping $Y_j$ to $f_jt, j=1,\ldots,r$ we obtain a representation of the Rees algebra:
$$A[It]\cong S/J,$$
where $J$ is the ideal of $A$ generated by the forms vanishing at $f_1,\ldots,f_r$.
If we set $\bideg x_i = (1,0)$ and $\bideg Y_j = (d_j,1)$, then $S$ is a bigraded algebra and the above isomorphism  is a bigraded isomorphism. 

For every $h=(\alpha_0,\ldots,\alpha_r)\in\NN^{r+1}$ put $S_h := A_{\alpha_0}Y_1^{\alpha_1}\ldots Y_r^{\alpha_r}$. Then $S =\oplus_{h\in\NN^{r+1}}S_h$ is an $\NN^{r+1}$-graded algebra. This $\NN^{r+1}$-grading is finer than the above bigrading because
$$S_{(u,v)} = \bigoplus_{\alpha_0+a_1\alpha_1+\ldots+ a_r\alpha_r=u\atop
\alpha_1+\ldots+\alpha_r=v}  S_{(\alpha_0,\alpha_1,\ldots,\alpha_r)},$$
for all $(u,v)\in\NN^2$. 
We order $\NN^{r+1}$ as follows: $(\alpha_0,\ldots,\alpha_r) < (\beta_0,\ldots,\beta_r)$ if the first non-zero component from the left side of 
$$\big(\alpha_0 + \sum_{j=1}^r\alpha_jd_j- \beta_0-\sum_{j=1}^r\beta_jd_j, \sum_{j=1}^r\alpha_j -\sum_{j=1}^r\beta_j,\alpha_0-\beta_0,\ldots,\alpha_r-\beta_r\big)$$
is negative. Then $<$ is a term order on $\NN^{r+1}$, i.e., $h < h'$ implies $h+g < h'+g$ for any $g \in \NN^{r+1}$. Note that this term order is different from that of [HeTU].

For every polynomial $f \in S$ let $f^*$ denote the initial term of $f$, i.e. $f^*:= f_h$ if $f =\sum_{h'\in \NN^{r+1}}f_{h'}$ and $h =\min \{h'|\ f_{h'}\neq 0\}$. Let $J^*$ denote the ideal of $S$ generated by the elements $f^*$, $f\in J$. Then $S/J^*$ is a bigraded algebra. This algebra has  a simpler structure than that of $S/J$. We can use $S/J^*$ to compute the Hilbert function $H_{A[It]}(u,v)$ by the following lemma.

\begin{Proposition} \label{approximation}
$H_{A[It]}(u,v) = H_{S/J^*}(u,v)$ for all $(u,v) \in \NN^2$.
\end{Proposition}

\begin{pf} Fix $(u,v) \in \NN^{r+1}$. Let
$$D  := \{(\alpha_0,\ldots,\alpha_r)\in \NN^{r+1}|\ \alpha_0+a_1\alpha_1+\ldots+ a_r\alpha_r=u,\
\alpha_1+\ldots+\alpha_r=v\}.$$
Then $S_{(u,v)} = \oplus_{h \in D}S_h$. Hence $S_{(u,v)} = 0$ if $D  = \emptyset$. If $D  \neq \emptyset$, we set
\begin{align*} 
h_m & := \min\{h|\ h \in D\},\\
h_M & := \max\{h|\ h \in D\}.
\end{align*}
By the definition of the order $<$ we have
$D = \{h \in \NN^{r+1}|\ h_m \le h \le h_M\}.$ 
For every $h \in \NN^{r+1}$ let 
\begin{align*}
F_h & := \oplus_{h' \ge h}S_{h'},\\
h^* & := \min\{h' \in \NN^{r+1}|\ h' > h\}.
\end{align*}
Then $J^*_h \cong J \cap F_h/J \cap F_{h^*}$. Moreover,
$F_{h_m} = \sum_{h \in D}S_h \oplus F_{(h_M)^*}= S_{(u,v)} \oplus F_{(h_M)^*}.$
This implies 
$J_{(u,v)} = J \cap F_{h_m}/J \cap F_{(h_M)^*}.$
Using the chain $J \cap F_{h_m} \subset J \cap F_{(h_m)^*} \subset \ldots \subset J \cap F_{h_M} \subset J \cap F_{(h_M)^*}$ we get
\begin{align*}
\dim_k J_{(u,v)} & = \sum_{h_m \le h \le h_M} \dim_k(J \cap F_h/J \cap F_{h^*})\\ & =\sum_{h \in D} \dim_kJ^*_h = \dim_k \bigoplus_{h \in D}J^*_h\\
& = \dim_k J^*_{(u,v)}.
\end{align*}
Hence $\dim_k (S/J)_{(u,v)} = \dim_k (S/J^*)_{(u,v)}$. So we get
$H_{A[It]}(u,v) = H_{S/J^*}(u,v)$.
\end{pf} 

\section{Rees algebras of homogeneous $d$-sequences}

Let $A=\oplus_{u\ge 0}A_u$ be a standard graded algebra over a field $k$. Let $f_1,\ldots,f_r$ be a sequence of homogeneous elements of $A$ and $I=(f_1,\ldots,f_r)$. 

We call $f_1,\ldots,f_r$ a $d$-{\it sequence} if the following conditions are satisfied:\par

(1) $f_i\notin (f_1,\ldots f_{i-1},f_{i+1},\ldots f_r),$\par

(2) $(f_1,\ldots f_i):f_{i+1}f_j=(f_1,\ldots,f_i):f_j$ for all $j \ge i+1$ and all $i\ge 0.$ \par

\noindent This notion was introduced by Huneke in [Hu], where one can find basic properties and abundant examples of $d$-sequences. In this section we will estimate the Hilbert polynomial of the bigraded Rees algebra $A[It]$ when $I$ is generated by a $d$-sequence of homogeneous elements of increasing degrees.  

By Proposition \ref{approximation} we only need to compute the Hilbert polynomial of $S/J^*$, where $S = A[Y_1,\ldots,Y_r]$ and $J^*$ is the initial ideal of the defining ideal of $A[It]$ with respect to the term order introduced there. 
The following result on the initial ideal $J^*$ is similar to [HeTU, Lemma 1.2] (which uses a different term order). It can be proved by the same proof. We leave the reader to check the details. 

\begin{Lemma} \label{ini} 
Let $f_1,\ldots,f_r$ be a homogeneous $d$-sequence with $\deg f_1\le\ldots\le\deg f_r$. Then
$$J^*=(I_1Y_1,\ldots,I_rY_r),$$
where $I_q=(f_1,\ldots,f_{q-1}):f_q$ for $q = 1,\ldots,r.$
\end{Lemma}

Using Lemma \ref{ini} we can compute the Hilbert polynomial of $S/J^*$ and get the following result on the Hilbert polynomial of $A[It]$.  

\begin{Theorem} \label{d-sequence}
Let $I$ be an ideal generated by a homogeneous $d$-sequence $f_1,\ldots,f_r$ with $\deg f_j= d_j$ and $d_1\le\ldots\le d_r$. Let $I_q=(f_1,\ldots,f_{q-1}):f_q$ for $q=1,\ldots,r$. Set 
$s := \dim A/I_1-1$ and $m := \max\{q|\ \dim A/I_q + q-2 = s\}$.
Then $\deg P_{A[It]}(u,v) = s$ and
$$e_i(A[It]) =  \sum_{q=1}^{\min\{m,s-i+1\}}(-1)^{s-q-i+1}e(A/I_q)\sum_{j_1+\ldots+j_q=s-q-i+1}d_1^{j_1}\ldots d_q^{j_q}$$
for $i = 0,\ldots,s$.
\end{Theorem}

\begin{pf} As observed above, we may replace $A[It]$ by $S/J^*$. We have
$$(S/J^*)_{(u,v)}=\bigoplus_{\alpha_0+d_1\alpha_1+\ldots+d_r\alpha_r=u\atop \alpha_1+\ldots+\alpha_r=v} (S/J^*)_{(\alpha_0,\alpha_1,\ldots,\alpha_r)}.$$
Every element of $J^*_{(\alpha_0,\alpha_1,\ldots,\alpha_r)}$ has the form $fY_1^{\alpha_1}\ldots Y_r^{\alpha_r}$ with $f \in (I_j)_{\alpha_0}$ for some $j=1,\ldots,r$ with $\alpha_j\not=0$. Since $I_1 \subset \ldots \subset I_r$ [Hu, Remarks (2)], this implies $$J^*_{(\alpha_0,\alpha_1,\ldots,\alpha_r)}=(I_{m(\alpha_1,\ldots,\alpha_r)})_{\alpha_0}Y_1^{\alpha_1}\ldots Y_r^{\alpha_r}$$
where 
$m(\alpha_1,\ldots,\alpha_r):=\max{\{q|\,\alpha_q\neq 0}\}$.  Hence
$$(S/J^*)_{(\alpha_0,\alpha_1,\ldots,\alpha_r)} \cong (A/I_{m(\alpha_1,\ldots,\alpha_r)})_{\alpha_0}.$$
So we obtain 
\begin{align*}
H_{S/J^*}(u,v) & =\sum_{\alpha_0+d_1\alpha_1+\ldots+d_r\alpha_r=u\atop \alpha_1+\ldots+\alpha_r=v}H_{A/I_{m(\alpha_1,\ldots,\alpha_r)}}(\alpha_0)\\
& =\sum_{\alpha_1+\ldots+\alpha_r=v}H_{A/I_{m(\alpha_1,\ldots,\alpha_r)}}(u-d_1\alpha_1-\ldots-d_r\alpha_r)\\
& =\sum_{q=1}^r \sum_{\alpha_1+\ldots+\alpha_q=v\atop \alpha_q\neq 0}H_{A/I_q}(u-d_1\alpha_1-\ldots-d_q\alpha_q)\\
& =\sum_{q=1}^r \big[\sum_{\alpha_1+\ldots+\alpha_q=v}H_{A/I_q}(u-d_1\alpha_1-\ldots-d_q\alpha_q)-\\
&\hspace{1.5cm} \sum_{\alpha_1+\ldots+\alpha_{q-1}=v}H_{A/I_q}(u-d_1\alpha_1-\ldots-d_{q-1}\alpha_{q-1})\big]
\end{align*}
		
By Theorem \ref{exist} there exist integers $u_0$ and $v_0$ such that $H_{S/J^*}(u,v) = P_{S/J^*}(u,v)$ for $u \ge dv + u_0$ and $v \ge v_0$. If $\alpha_1+\ldots+\alpha_j = v$, then $u-d_1\alpha_1-\ldots-d_j\alpha_j \ge  u-dv$ so that we may choose $u_0$ such that
$$H_{A/I_q}(u-d_1\alpha_1-\ldots-d_j\alpha_j) = P_{A/I_q}(u-d_1\alpha_1-\ldots-d_j\alpha_j)$$
for $u \ge dv+u_0$. Put 
$H_{q,j}(u,v) := \displaystyle \sum_{\alpha_1+\ldots+\alpha_j=v}P_{A/I_q}(u-d_1\alpha_1-\ldots-d_j\alpha_j).$ Then
$$P_{S/J^*}(u,v) = \sum_{q=1}^r \big[H_{q,q}(u,v)- H_{q,q-1}(u,v)\big].$$
For $q = 1,\ldots,r$ let $s_q := \dim A/I_q-1$.  By Lemma \ref{observation}, $H_{q,j}(u,v)$ is a polynomial of degree $s_q+j-1$.
Therefore,
$$\deg P_{S/J^*}(u,v) \le \max\{s_q+q-1|\ q = 1,\ldots,r\}.$$ 
Since $I_q \supseteq (I_{q-1},f_{q-1})$ and $f_q$ is a non-zerodivisor modulo $I_{q-1}$, we have
$$s_q  < \dim A/(I_{q-1},f_{q-1})  = \dim A/I_{q-1}-1= s_{q-1}$$
for $q \ge 2$.  From this it follows that
$s_q+q-1 \le s_1 =s$. Hence $\deg P_{S/J^*}(u,v) \le s$.\par

Since $s = s_1 > \ldots > s_r$ and $m = \max\{q|\ s_q+q-1 = s\}$, we have $s_q+q-1 = s$ if $q \le m$ and $s_q + q-1 < s$ if $q > m$. Therefore, $P_{S/J^*}(u,v)$ and $\sum_{q=1}^{m}H_{q,q}(u,v)$ share the same terms of degree $s$. Using Lemma \ref{observation} we can compute these terms and we may write
\begin{align*} 
P_{S/J^*}(u,v)  & = \sum_{q=1}^m\sum_{i=0}^{s_q}\frac{(-1)^{s_q-i}e(A/I_q)\displaystyle \sum_{j_1+\ldots+j_q=s_q-i}d_1^{j_1}\ldots d_q^{j_q}}{i!(s-i)!}u^iv^{s-i} +\\
& \hspace{1.5cm} +\text{\rm lower-degree terms}.
\end{align*}
For $q = 1,\ldots,m$ we have $s_q = s-q+1$. Therefore, $P(u,v)$ can be rewritten as 
\begin{align*}
P_{S/J^*}(u,v)  & = \sum_{i=0}^{s}\frac{\displaystyle\sum_{q=1}^{\min\{m,s-i+1\}}(-1)^{s-q-i+1}e(A/I_q)\sum_{j_1+\ldots+j_q=s-q-i+1}d_1^{j_1}\ldots d_q^{j_q}}{i!(s-i)!}u^iv^{s-i} +\\
& \hspace{1.5cm} + \text{\rm lower-degree terms}.
\end{align*}
The coefficient of $u^s$ in $P(u,v)$ is equal to $e(A/I_1)$. So we can conclude that
$\deg P_{S/J^*}(u,v)  = s$ and 
$$e_i(S/J^*) = \sum_{q=1}^{\min\{m,s-i+1\}} (-1)^{s-q-i+1} e(A/I_q)\sum_{j_1+\ldots+j_q=s-q-i+1}d_1^{j_1}\ldots d_q^{j_q}.$$
\end{pf}

\noindent{\it Remark.}  The above method was already used in [RaV] (see also [Ho]) to establish the following formula for $e_i(\mm|I)$ when $I$ is generated by a $d$-sequence: 
$$e_i(\mm|I)= \cases 0 & \text {if $0\le i\le s-m,$}\\
e(R/{I_{s-i+1}}) & \text {if $s-m < i \le s.$}
\endcases$$
This displays a completely different behavior than that of $e_i(A[It])$.  

Now we will apply Theorem \ref{d-sequence} to some special class of $d$-sequences.  

\begin{Corollary} \label{regular} Let $f_1,\ldots f_r$ be a homogeneous regular sequence with $\deg f_1=d_1 \le\ldots\le \deg f_r = d_r$ and $I = (f_1,\ldots,f_r)$. Set $s = \dim A-1$. Then $\deg P_{A[It]}(u,v)  = s$ and 
$$e_i(A[It]) = \sum_{q=1}^{\min\{r,s-i+1\}} (-1)^{s-q-i+1} e(A)\sum_{j_1+\ldots+j_q=s-q-i+1}d_1^{j_1+1}\ldots d_{q-1}^{j_{q-1}+1}d_q^{j_q}$$
for $i = 0,\ldots,s$.
\end{Corollary}

\begin{pf} Since $I_q = (f_1,\ldots,f_{q-1}):f_q = (f_1,\ldots,f_{q-1})$, we have 
$\dim A/I_q  = \dim A -q+1$ and $e(A/I_q) = d_1\ldots d_{q-1}e(A)$, $ q = 1,\ldots,r$. It follows that $s = \dim A/I_1-1$ and $r = \max\{q|\ \dim A/I_q +q -2 = s\}$. Hence the conclusion follows from Theorem \ref{d-sequence}.
\end{pf}

If $A$ is a polynomial ring and $I$ a completion ideal in $A$,  the Hilbert function of $A/I_q$ can be expressed in terms of $d_1,\ldots,d_q$, $q = 1,\ldots,r$. Therefore, we can use the formula for $H_{S/J^*}(u,v)$ in the proof of Theorem \ref{d-sequence} to give an explicit formula for $H_{A[It]}(u,v)$ in terms of $d_1,\ldots,d_r$. We will demonstrate this method by the following example which covers a key result of Geramita, Gimigliano and Harbourne in [GGH].  

\begin{Example} \label{GGH}
{\rm Let $f_1,f_2$ be a homogeneous regular sequence in $A = k[x_1,x_2,x_3]$ with $\deg f_1 = d_1 \le \deg f_2 = d_2$ and $I = (f_1,f_2)$. The rational surface $V$ obtained by  blowing up the intersection of the curves $f_1$ and $f_2$ was the main focus of study in [GGH], where it was assumed that $f_1$ and $f_2$ meet transversally in $d_1d_2$ points.  The key result there is an explicit formula for the cohomological function $h^0(V,D_{u'v})$ in terms of $d_1,d_2$, where $D_{u',v}$ denotes the divisor class of $V$ associated to $(I^v)_{u'+d_2v}$ [GGH, Proposition III.1(a)]. This formula can be easily derived by our method. 
We first note that
$$h^0(V,D_{u'v}) = \dim_k (I^v)_{u'+d_2v} = H_{A[It]}(u'+d_2v,v).$$
Since $I_1 = 0$, $I_2 = (f_1)$ and $H_A(u) = {u+2 \choose 2}$, $H_{A/(f_1)}(u) = {u+2\choose 2}-{u-d_1+2 \choose 2}$, we have
\begin{align*}
H_{A[It]}(u,v) & = H_A(u-d_1v) +\sum_{j=0}^vH_{A/(f_1)}(u-d_1j-d_2v+d_2j)- H_{A/(f_1)}(u-d_1v)\\
& =  \sum_{j=0}^v\left[{u-d_1j-d_2v+d_2j +2\choose 2}- {u-d_1j-d_2v+d_2j -d_1+2\choose 2}\right]+\\
&\quad + {u-d_1v-d_1+2 \choose 2}
\end{align*}
for $u > d_2v$. Putting $u = u'+d_2v$ and $\delta = d_2-d_1$ we get
\begin{align*}
H_{A[It]}(u'+d_2v,v) & = \sum_{j=0}^{v-1}\left[{u'+\delta j+2 \choose 2}-{u'+\delta j - d_1+2\choose 2}\right] +{u'+\delta v +2 \choose 2}.
\end{align*}
This is precisely the formula of [GGH, Proposition III.1(a)] for $h^0(V,D_{u'v})$. Moreover, using Corollary \ref{regular} we can easily compute the leading coefficients of the above polynomial:
$$e_0(A[It]) = - d_1d_2,\ e_1(A[It]) = 0, \ e_2(A[It])  = 1.$$ 
By Corollary \ref{embedded}, this implies the following formula for the degree of the embedded variety $V_{c,e}$ of $V$:
$$\deg V_{c,e} = c^2 - d_1d_2e.$$}
\end{Example}

\begin{Corollary} \label{minor} Let $A = k[X]$, where $X$ is a $(r-1)\times r$ matrix of indeterminates. Let $I$ be the ideal of the maximal minors of $X$ in $A$. Set $s = (r-1)\times r-1$. Then $\deg P_{A[It]}(u,v)  = s$ and 
$$e_i(A[It]) = \sum_{q=1}^{\min\{r,s-i+1\}} (-1)^{s-q-i+1}{r -1\choose q-1} {s-i \choose q-1} r^{s-q-i+1}$$
for $i = 0,\ldots,s$.
\end{Corollary}

\begin{pf} Let $f_q$ be the determinant of the submatrix of $X$ obtained by deleting the $(r-q+1)$th column, $q = 1,\ldots,r$. Then $\deg f_1 = \ldots = \deg f_r = r$ and $f_1,\ldots,f_r$ is a $d$-sequence by [Hu, Example 1.1]. Let $I_q = (f_1,\ldots,f_{q-1}):f_q$. Then $I_q$ is the ideal generated by the maximal minors of the matrix of the first $r-q+1$ columns of $X$. Hence $\dim R/I_q = (r-1)r-q+1$ and $e(A/I_q) = {r -1\choose q-1}$ [HeT, Theorem 3.5]. It follows that $s = \dim R/I_1-1$ and $r = \max\{q|\ \dim A/I_q + q -2 = s\}$. Since 
$$\sum_{j_1+\ldots+j_q=s-q-i+1}r^{j_1}\ldots r^{j_q} = {s-i \choose q-1} r^{s-q-i+1},$$
the conclusion follows from Theorem \ref{d-sequence}.
\end{pf}

\section*{References}

\noindent [B] P.B. Bhattacharya, The Hilbert function of two ideals, Proc. Cambridge Phil. Soc. 53 (1957), 568-575.\par
\noindent [CHTV] A. Conca, J. Herzog, N.V. Trung, G. Valla, Diagonal subalgebras of bigraded algebras and embeddings of blow-ups of projective spaces, Amer. J. Math. 119 (1997), 859-901.\par
\noindent [CH] S. D. Cutkosky and J. Herzog, Cohen-Macaulay coordinate rings of blowup schemes, Comment. Math. Helvetici 72 (1997), 605-617. \par
\noindent [E] D. Eisenbud, Commutative algebra with a view toward algebraic geometry,  Springer-Verlag, 1995.\par
\noindent [GG] A.V. Geramita and A. Gimigliano, Generators for the defining ideal of certain rational surfaces, Duke Math. J. 62 (1991), 61-83. \par
\noindent [GGH] A. V. Geramita, A. Gimigliano, and B. Harbourne,
Projectively normal but superabundant embeddings of rational
surfaces in projective space, J.~Algebra 169 (1994), 791--804.\par
\noindent [Gi]  A. Gimigliano, On Veronesean surfaces, Proc. Kon. Nederl. Akad. Wetensch. Ser. A 92 (1989), 71-85. \par
\noindent [GL] A. Gimigliano and A. Lorenzini, On the ideal of Veronesean surfaces, Canad. J. Math. 45 (4) (1993), 758-777.\par
\noindent [HeT] J. Herzog and N.V. Trung, Gr\"obner bases and multiplicity of determinantal and Pfaffian ideals, Adv. in Math. 96 (1992), 1-37.\par
\noindent [HeTU] J. Herzog, N.V. Trung and B. Ulrich,  On the multiplicity of Rees algebras and associated graded rings of $d$-sequences,  J. Pure Appl. Algebra~80 (1992), 273-297. \par
\noindent [Ho] N. D. Hoang, On mixed multiplicities of homogeneous ideals, Beitr\"age Algebra Geom. 42 (2001), 463--473. \par
\noindent [Hu] C. Huneke, The theory of d-sequence and powers of ideals, Adv. in Math. 46 (1982), 249-297.\par
\noindent [KMV] D. Katz, S. Mandal and J. Verma, Hilbert function of bigraded algebras, in:  A. Simis, N.V. Trung and G. Valla  (eds.), Commutative Algebra (ICTP, Trieste, 1992), 291-302, World Scientific, 1994.\par
\noindent [KV] D. Katz and J.K. Verma, Extended Rees algebras and mixed multiplicities, Math. Z. 202 (1989), 111-128.\par
\noindent [Ma] H. Matsumura, Commutative ring theory, Cambridge Studies in Advanced Math. 8,  Cambridge, 1989.\par 
\noindent [RaS] K.N. Raghavan and A. Simis, Multiplicities of blow-ups of homogeneous ideals generated by quadratic sequences, J. Algebra 175 (1995), 537-567.\par
\noindent [RaV] K.N. Raghavan and J.K. Verma, Mixed Hilbert coefficients of homogeneous $d$-sequences and quadratic sequences, J. Algebra 195 (1997), 211-232. \par
\noindent [Re] D. Rees, Generalizations of reductions and mixed multiplicities, J. London Math. Soc. 29 (1984), 423-432.\par
\noindent [Ro1] P. Roberts, Multiplicities and Chern classes in local algebra, Cambrdige Tracts in Math. 133, Cambridge,
1998. \par 
\noindent [Ro2] P. Roberts, Recent developments on Serre's multiplicity conjectures: Gabber's proof of the nonnegativity conjecture, L'Enseignment Math\'ematique 44 (1998), 305-324.\par
\noindent [Ro3] P. Roberts,  Intersection multiplicities and Hilbert polynomials,
Michigan Math. J. 48 (2000), 517-530. \par
\noindent [STV] A. Simis, N.V. Trung and G. Valla, The diagonal subalgebras of a blow-up ring, J. Pure Appl. Algebra 125 (1998), 305-328. \par
\noindent [St] R. Stanley, Hilbert functions of graded algebras, Adv. in Math. 28 (1978), 57-83.\par
\noindent [SV] J. St\"uckrad and W. Vogel, Buchsbaum rings and applications, VEB Deutscher Verlag der Wisssenschaften, Berlin, 1986.\par
\noindent [Te] B. Teissier, Cycles \'evanescents, sections planes, et conditions de Whitney, Singularit\'es \`a Carg\`ese 1972, Ast\`erisque 7-8 (1973), 285-362. \par
\noindent [Tr1] N.V. Trung,  Filter-regular sequences and multiplicity of blow-up rings of ideals of the principal class,  J. Math. Kyoto Univ.~33 (1993), 665-683.\par
\noindent [Tr2] N.V. Trung, Positivity of mixed multiplicities, Math. Ann. 319 (2001), 33-63.\par
\noindent [Va] G. Valla, Certain graded algebras are always Cohen-Macaulay, J. Algebra 42 (1976), 537-548.\par
\noindent [Ve1] J.K. Verma, Rees algebras and mixed multiplicities, Proc. Amer. Math. Soc. 104 (1988), 1036-1044. \par
\noindent [Ve2] J.K. Verma, Multigraded Rees algebras and mixed multiplicities, J. Pure Appl. Algebra 77 (1992), 219-228. \par
\noindent [W] B.L. van der Waerden, On Hilbert's function, series of composition of ideals and a generalization of the theorem of Bezout, Proc. Kon. Nederl. Akad. Wetensch. Amsterdam 31 (1928), 749-770.
\end{document}